\documentclass[reqno,11pt]{amsart}
\usepackage{latexsym,amssymb,amsthm,amsmath}

\theoremstyle{plain}
\newtheorem{theorem}{Theorem}[section]
\newtheorem{proposition}{Proposition}[section]
\newtheorem{lemma}{Lemma}[section]
\newtheorem{corollary}{Corollary}[section]

\theoremstyle{definition}

\theoremstyle{remark}
\newtheorem{remark}{Remark}[section]
\newtheorem{example}{Example}[section]

\numberwithin{equation}{section}
\def\({\left( }
\def\){\right)}

\begin{document}
\setcounter{page}{1}

\title[2-ranks of compact Lie groups]{The 2-ranks of connected compact Lie groups}

\author[B.-Y. Chen]{ Bang-Yen Chen}

\address{Department of Mathematics,
	Michigan State University, 619 Red Cedar Road, East Lansing, Michigan
48824--1027, U.S.A.}

\email{bychen@math.msu.edu}

\begin{abstract} The 2-rank of a compact Lie group $G$ is the maximal possible rank of the elementary 2-subgroup ${\mathbb Z}_{2}\times\cdots {\mathbb Z}_{2}$ of $G$. The study of 2-ranks (and $p$-rank for any prime $p$) of compact Lie groups was initiated in 1953 by A. Borel and J.-P. Serre \cite{BS}. Since then the 2-ranks of compact Lie groups have been investigated by many mathematician. The 2-ranks of compact Lie groups  relate closely with several important areas in mathematics.

In this article, we survey important results concerning 2-ranks of compact Lie groups. In particular, we present the complete determination of 2-ranks of compact connected simple Lie groups $G$ via the maximal antipodal sets $A_{2}G$ of $G$ introduced  in \cite{CN2,CN3}.
\end{abstract}



\maketitle

\section{Introduction} 

For a prime number $p$,  a $p$-group is a periodic group in which each element has a power of $p$ as its order. 
A Sylow $p$-subgroup of a finite group $G$ is a maximal $p$-subgroup of $G$, i.e., a subgroup of $G$ that is a $p$-group, and that is not a proper subgroup of any other $p$-subgroup of $G$.

 In the field of finite group theory, the Sylow theorems are a collection of theorems named after  Ludwig Sylow (1832-1918) that give detailed information about the number of subgroups of fixed order that a given finite group contains (cf. \cite{S}). The Sylow theorems assert a partial converse to Lagrange's theorem that for any finite group $G$ the order of every subgroup of $G$ divides the order of $G$.

 The Sylow theorems form a fundamental part of finite group theory and have very important applications in the classification of finite simple groups (cf.  \cite{A,M}). Further, the problem of finding a Sylow subgroup of a given group becomes an important problem in computational group theory. In permutation groups, it has been proven in  that a Sylow $p$-subgroup and its normalizer can be found in polynomial time of the input. These algorithms  are now becoming practical as the constructive recognition of finite simple groups becomes a reality (cf. \cite{Ka,KT}). In particular, versions of this algorithm are now used in the Magma computer algebra system.

 A Sylow $p$-subgroup in an infinite group to be a $p$-subgroup  that is maximal for inclusion among all $p$-subgroups in the group. Such subgroups exist by Zorn's lemma. There is an analogue of the Sylow theorems for infinite groups.

 A. Borel and J.-P. Serre considered in \cite{BS} a topological group $H$ which has a sequence of normal subgroups $$e=H_{0}\subset H_{1}\subset \cdots\subset H_{k}=H$$  such that every factor group $H_{i}/H_{i-1}$ is either a finite cyclic group or a one-dimensional torus, and they prove that if such a group $H$ is a subgroup of a compact Lie group $G$, it is contained in the normalizer $N$ of a maximal torus $T$ in $G$. It follows, in particular, that every abelian subgroup of $G$ is contained in some $N$. 
 
Borel and Serre defined in \cite{BS} the $p$-rank $r_{p}$ ($p$ prime) of a compact Lie group $G$ as the largest integer $h$ such that $G$ contains the direct product of $h$ cyclic groups of order $p$. The $p$-ranks (in particular, the 2-ranks) of compact Lie groups relate closely with several important areas in mathematics.

It is clear that the $p$-rank of $G$ is at least equal to the rank of $G$ (i.e., the dimension of a maximal torus $T$), but since the Weyl group $N/T$ ($N$ the normalizer of $T$) of $G$ is finite, it follows from the above result that the $p$-rank of $G$ is always finite. 
 By applying the method of spectral sequences of fibre bundles, they 
 proved that if the $p$-rank of a connected compact Lie group $G$ is greater than the rank of $G$, then $G$ has $p$-torsion. 
 
 In particular, for $2$-rank, A. Borel and J.-P. Serre established the following two results:

\begin{itemize}
\item[(i)] The usual rank $rank(G)\leq r_{2}G\leq 2 ( rank(G))$; and

\item[(ii)] $G$ has (topological) 2-torsion if $rank(G)<r_{2}G$.
\end{itemize}

In \cite{BS}, they are able to determine the 2-rank of the simply-connected simple Lie groups $SO(n), Sp(n), U(n), G_{2},F_{4}.$ They showed that the exceptional Lie groups $G_{2}$, $F_{4}$ and $E_{8}$ have 2-torsion. On the other hand, it was also mentioned in \cite[page 139]{BS} that they are unable to determine the 2-rank for the exceptional simple Lie groups $E_{6}$ and $E_{7}$. Since then the 2-ranks of compact Lie groups have been investigated by many mathematician. 
 
In this article, we survey important results concerning the 2-rank of connected compact Lie groups. In particular, we provide the complete determination of the 2-ranks of compact connected simple Lie groups $G$ via the maximal antipodal sets $A_{2}G$ and the  2-number $\#_{2}G$ of $G$, which were originally introduced by B.-Y. Chen and T. Nagano in \cite{CN2,CN3}. The maximal antipodal sets of compact simple Lie groups were determined in \cite{CN3} via the $(M_{+},M_{-})$-method discovered in \cite{CN1}.

 The main references of this article are \cite{BS,c,CN1,CN2,CN3,Q1,Q2}.

 \section{Basic notations and definitions}  
 
 In this article, we assume a good knowledge of \cite{H} on Lie groups and symmetric spaces on the
part of the readers.

 \subsection{Symmetric spaces}
 
 An isometry $s$ of a Riemannian manifold $M$ is said to be involutive if $s^{2}=id_{M}$. A Riemannian manifold $M$ is called a symmetric space if for each point $x\in M$ there is an involutive isometry $s_{x}$ of $M$ such that $x$ is an isolated fixed point of $s_{x}$. The involutive isometry $s_{x}$ is called the symmetry at $x$.
 
 For a symmetric space $M$, we denote by $G_{M}$ the identity component of the isometric group $I(M,g)$ of $M$. We denote by $K_{M}$ the isotropy subgroup at a point $o$ so that $M=G_{M}/K_{M}$.
 
  For a compact symmetric space $M$, let  $M^{*}$ denote  the bottom space (the adjoint space in \cite{H}) of the
space $M$. Let $M/{\mathbb Z}\mu$ denote the space of which $M$ is a covering space with the covering
transformation group ${\mathbb Z}\mu$, the cyclic group of order $\mu$ (if there is no ambiguity).

 \subsection{Antipodal set and the 2-number} 

The notions of maximal antipodal sets and  2-number were introduced by  Chen and  Nagano in \cite{CN2,CN3}.

For a compact symmetric space $M$, the 2-number, denoted by $\#_2M$, is defined as
the maximal possible cardinality $\#_2A_2M$ of a subset $A_2M$ of $M$ such that the point symmetry $s_x$ fixes every point of $A_2M$ for every $x\in A_2M$. It is known that a subset $A_{2}$ of a connected symmetric space $M$ is antipodal if and only if any pair of points in $A_{2}$ are antipodal points on some closed geodesic in $M$.

The 2-number $\#_2M$ is finite. The definition is equivalent to saying that $\#_2M$ is a maximal possible cardinality $\#A_2M$ of a subset $A_2M$ of $M$ such that, for every pair of points $x$ and $y$ of $A_2M$, there is a closed geodesic of $M$ on which $x$ and $y$ are antipodal to each other. Thus the invariant can also be defined on Riemannian manifolds.

The invariant, $\#_2M$, has certain bearings on the topology of $M$; for instance, $\#_2M$ equals $\mathcal X(M)$, the Euler number of $M$, if $M$ is a semisimple hermitian symmetric space.  And in general one has the inequality $\#_2M\geq \mathcal X(M)$ for any compact connected symmetric space $M$. Furthermore, M. Takeuchi proved in \cite{T} that $\#_{2}M=\dim H(M,\mathbb Z_{2})$ for any symmetric $R$-space, where $H(M,\mathbb Z_{2})$ is the homology group of $M$ with coefficients in $\mathbb Z_{2}$. This formula is actually correct for every
space we could check (see \cite{CN3}).

\subsection{Compact Lie groups}
If $G$ is a connected compact Lie group, then by assigning $s_{x}(y)=xy^{-1}x$ to every point $x\in G$, we have $s^{2}_{x}=id_{G}$ to each point $x$. Therefore $G$ becomes a compact symmetric space with respect to a bi-invariant Riemannian metric. 
 
  The automorphism group, denoted by $Aut(G)$, of a connected compact Lie group $G$ contains the left translation group, and hence $Aut(G)$ is transitive on $G$.

Any compact connected Lie group $G$ contains a maximal torus $T$. The solutions of
$t^{2}=1$ in $T$ form a 2-subgroup $A_{2}T\subset G$ of the same rank as $G$. However, it is not
necessarily maximal. This fact shows  $rank(G)\leq r_{2}G$.

In this article, the standard notations for the Lie groups such as $G_{2}, F_{4},E_{6},\cdots $ denote the simply-connected ones.

 \subsection{$(M_+,M_-)$-method} Here we provide a  brief introduction of the $(M_+,M_-)$-method for compact symmetric spaces, introduced by B.-Y.  Chen and  T. Nagano in \cite{CN1,CN3}. This method plays the key roles in our determination of 2-numbers of compact symmetric spaces; in particular, of 2-ranks of compact Lie groups.

 Let $o$ be a point of a symmetric space $M$. We call a connected component of the fixed point set $F(s_{o},M)$ of the symmetry $s_{o}$ in $M$ a {\it polar} of $o$. We denote it by $M_{+}$ or $M_{+}(p)$ if $M_{+}$ contains a point $p$.  When a polar consists of a single point, we call it a {\it pole}. 
 
 We call a connected component of the fixed point set $F(s_{p}\circ s_{o},M)$ of $s_{p}\circ s_{o}$ through $p$ the {\it meridian} of $M_{+}(p)$ in $M$ and denote it by $M_{-}(p)$ or simply by $M_{-}$.  Notice that $M_{-}(p)$ has the same rank as $M$. Moreover, we have $$\dim M_{+}(p)+\dim M_{-}(p)=\dim M.$$

Polars and meridians are totally geodesic submanifolds of a symmetric space $M$; they are thus symmetric spaces. And they have been determined for every compact connected irreducible Riemannian symmetric space (see  \cite{c,CN1,Ki,Na,Na2}). One of the most important properties of these totally geodesic submanifolds is that $M$ is determined by any pair of $(M_{+}(p),M_{-}(p))$ completely. 

 Let $M$ be a compact connected Riemannian symmetric space and $o$ be a point in $M$. If there exists a pole $p$ of $o\in M$, then we call the set consisting of the midpoints of the geodesic segments from $o$ to $p$ the {\it centrosome} and denote it by $C(o,p)$. Each connected component of the centrosome of $\{o,p\}$ is a totally geodesic submanifold of $M$.

The following result from \cite{CN3} characterizes poles in compact symmetric spaces.

\begin{proposition} \label{P:2.1}
 The following six conditions are equivalent to each other for two distinct points $o,p$ of a connected compact symmetric space $M=G_{M}/K_{G}$.

\begin{itemize} 

\item[(i)] $p$ is a pole of $o\in M$; 
  
\item[(ii)] $s_{p}=s_{o}$; 
  
\item[(iii)] $\{p\}$  is a polar of $o\in M$; 
  
\item[(iv)] there is a double covering totally geodesic immersion $\pi=\pi_{\{o,p\}} :M\to M''$ with $\pi(p)=\pi(o)$;

\item[(v)] $p$  is a point in the orbit $F(\sigma,G_{M})(o)$  of the group $F(\sigma,G_{M})$ through $o$, where $\sigma = {\rm ad}(s_{o})$;

\item[(vi)] the isotropy subgroup of $SG_{M}$ at $p$ is that, $SK_{G}$ $($of $SG_{M}$ at $o)$, where $SG_{M}$ is the group generated by $G_{M}$ and the symmetries; $SG_{M}/G_{M}$ is a group of order $\leq 2$.
\end{itemize}  
\end{proposition}
  
 For a compact symmetric space $M$, the {\it Cartan quadratic morphism} $$Q=Q_{o}: M\to G_{M}$$ carries a point $x\in M$ into $s_{x} s_{o}\in G_{M}$. The Cartan  quadratic morphism is a $G_{M}$-equivariant morphism which is an immersion.
  
 We have the following result  for centrosomes \cite{CN3}.
  
\begin{proposition} \label{P:2.2}   The following five conditions are equivalent to each other for two distinct points $o,q$ of a connected compact symmetric space $M$.

\begin{itemize} 

\item[(i)] $s_{o}s_{q}=s_{q}s_{o}$; 
  
\item[(ii)] $Q(q)^{2}=1_{G_{M}}$, where $Q=Q_{o}$ is Cartan  quadratic morphism; 
  
\item[(iii)] either $s_{o}$ fixes $q$ or $q$ is a point in the controsome $C(o,p)$ for some pole $p$ of $o$; 
  
\item[(iv)] either $s_{o}(q)=q$ or $s_{o}(q)=\gamma (q)$ for the covering transformation $\gamma$ for some pole $p=\gamma(o)$ of $o$;

\item[(v)] either $s_{o}(q)=q$ or there is a double covering morphism $\pi:M\to M''$ such that $s_{o''}$ fixes $q''$, where $o''=\pi(o)$ and $q''=\pi(q)$.
\end{itemize}  
\end{proposition}

These two propositions plays important roles for the study of maximal antipodal sets and 2-numbers of compact symmetric spaces and of compact Lie groups.

It was known that the polars, meridians and centrosomes play very important roles in the study of compact symmetric spaces as well as of compact Lie groups (cf. \cite{c,c1,c2000,CN1,CN3}).

\section{Two basic results on 2-ranks}

 The relationship between the 2-rank and the 2-number of a compact Lie group $G$ was established in the following theorem of \cite{CN3}.
 
 \begin{theorem}\label{T:3.1} Let $G$ be a connected compact Lie group. Then we have
 \begin{equation} \#_{2}G=2^{r_{2}G}.\end{equation}
 \end{theorem}
 This theorem can be proved as follows: Since the automorphism group $Aut(G)$ is transitive on $G$, in finding the 2-number $\#_{2}G$ we may therefore assume that a maximal antipodal set $A_{2}$ in $G$ contains the unit element 1. First observe
that the fixed point set $F(s_{1}, G)$ consists of the members of order 1 or 2. Thus we
have $s_{x}(y) = y$ for any two members $x, y$ of $A_{2}$ if and only if $xy = yx$ holds. Also observe
that this implies that the maximal $A_{2}$ is a subgroup. Moreover, $A_{2}$ is an elementary
abelian 2-subgroup isomorphic to $({\mathbb Z}_{2})^{t}={\mathbb Z}_{2}\times\cdots {\mathbb Z}_{2}$ ($t$ copies of $\mathbb Z_{2}$) for some positive integer $t$, a 2-subgroup for short. The largest possible value of
$t$ is by definition the 2-rank of the Lie group $G$. Consequently, $\#_{2}G$ is
a power of 2 for the compact Lie group $G$.
 
 For products of compact Lie groups, we have the following.
 
\begin{theorem}\label{T:3.2} Let $G_{1}$ and $G_{2}$ be  connected compact Lie groups. Then 
 \begin{equation}\label{3.2} \#_{2}(G_{1}\times G_{2})=2^{r_{2}G_{1}+r_{2}G_{1}}.\end{equation}
 \end{theorem}
 
This follows from the fact that the symmetry at $(x,y)\in G_{1}\times G_{2}$ carries a point $(u,v)$ into $(s_{x}u, s_{y}v)$. Thus one has $F(s_{(x,y)},G_{1}\times G_{2})=F(s_{x},G_{1})\times F(s_{y}, G_{2})$. Consequently, we have \eqref{3.2}.

  \section{Relations between 2-ranks and topology}

The first relationship between the 2-number and the topological 2-torsion of a compact Lie group was discovered by A. Borel and J.-P. Serre in \cite{BS}, in which they applied the method of spectral sequences of fiber bundles.

 \begin{theorem}\label{T:4.1} Let $G$ be a connected compact Lie group. If \begin{equation} rank (G)<2^{r_{2}G},\end{equation} then $G$ has topological 2-torsion.
 \end{theorem}

By  using the spectral sequence connecting $H^{2}(G,{\mathbb Z}_{2})$ to $H^{2}(G,{\mathbb Z})/{\rm Tors} \otimes {\mathbb Z}_{2}$, whose differential is the successive Bockstein operators, A. Borel also obtained  the following two results  in \cite{B2}.

 \begin{theorem}\label{T:4.2} The cohomology of a connected compact Lie group $G$ has no topological
2-torsion if and only if every antipodal subgroup $A_{2}G$ is contained in some torus $T$ in $G$. \end{theorem}

 \begin{theorem}\label{T:4.3}  If a compact Lie group $G$ is simply-connected and
it does have a topological 2-torsion, then $G$ contains an antipodal group  of rank 3 which no
torus of $G$ contains. \end{theorem}

By applying the classification of simply-connected simple Lie groups and on case-by-case consideration,
 A. Borel also proved  in \cite{B3} the following two results which  relate the 2-subgroups and 2-torsions for compact Lie groups.

\begin{theorem}\label{T:4.4} Let $G$ be a simply-connected simple Lie group. Then $G$ has no topological $2$-torsion  if $2$ does not divide the coefficients in the expression of the highest root as a linear combination of simple roots. 
\end{theorem}

\begin{theorem}\label{T:4.5} A connected compact simple Lie group $G$ has no topological 2-torsion if and only if every 2-subgroup  is contained in a torus. 
\end{theorem}

The following simple relationship between the 2-rank and the Euler characteristic was discovered by Chen and Nagano in \cite{CN3}.

 \begin{theorem}\label{T:4.6} Let $G$ be a connected compact Lie group. Then we have \begin{equation} 2^{r_{2}G}\geq {\mathcal X} (G),\end{equation} where $\mathcal X(G)$ denotes the Euler characteristic of $G$.
  \end{theorem}

Theorem \ref{T:4.6} was proved by applying the $(M_{+},M_{-})$-method together with some results of H. Hopf and H. Samelson \cite{HS}.

The following result was obtained by M. Takeuchi in \cite{T} who applied a result of Chen-Nagano in \cite{CN3} as well as an earlier result of Takeuchi \cite{T1965}.

 \begin{theorem}\label{T:4.7} If $M$ is a symmetric $R$-space, then
  \begin{equation}\label{4.3}\#_{2}M=\dim H(M,\mathbb Z_{2}),\end{equation} 
  where $H(M,\mathbb Z_{2})$ is the homology group of $M$ with coefficients in $\mathbb Z_{2}$. 
  
  In particular, \eqref{4.3} holds for the classical simple Lie groups $SO(m), U(m)$ and $Sp(m)$.
  \end{theorem}

 Let $F$ be either a field or the rational integer ring $\mathbb Z$. Let $$A=\sum_{i\geq 0} A_{i}$$ be a graded commutative $F$-algebra in sense of Milnor-Moore \cite{MM}. If $A$ is connected, then it admits a unique augmentation $\varepsilon :A\to F$.   Put  $\bar A={\rm Ker}\, \varepsilon$. The $\bar A$ is called the augmentation ideal of $A$.

  A sequence of elements $\{x_{1},\ldots,x_{n}\in \bar A\}$ in  the augmentation ideal $\bar A$ is called a {\it simple system of generators} if $\{x_{1}^{\epsilon_{1}}\cdots x_{n}^{\epsilon_{n}}: \epsilon_{i}=0 \hbox{ or }1\}$ is a module base of $A$.

 Let $G$ be a compact connected Lie group. Denote by $s(G)$ the number of generators of a simple system of the $\mathbb Z_{2}$-cohomology $H^{*}(G,{\mathbb Z}_{2})$ of  $G$. 

A. Kono discovered in \cite{K} the following relationship between $s(G)$,  $r_{2}G$ and the $\mathbb Z_{2}$-cohomology $H^{*}(G,{\mathbb Z}_{2})$.

 \begin{theorem}\label{T:4.8} Let $G$ be a connected compact Lie group. Then
  the following three conditions are equivalent: 
  
  \begin{itemize}
  \item[{\rm (1)}] $s(G)\leq r_{2}G$; 
  
  \item[{\rm (2)}] $s(G)= r_{2}G$; 
  
  \item[{\rm (3)}] $H^{*}(G,{\mathbb Z}_{2})$ is generated by universally transgressive elements. 
\end{itemize}
\end{theorem}

To prove Theorem \ref{T:4.8}, A. Kono used May's spectral sequence \cite{May},  Eilenberg-Moore's spectral sequence \cite{EM} as well as Quillen's result in \cite{Q1}. 

In \cite{K}, Kano also described properties of compact Lie groups satisfying condition (3) in Theorem \ref{T:4.8} and gives some applications.

\section{Covering maps and 2-ranks}

In this section, we present some very simple relationships between covering maps and 2-numbers for connected compact Lie groups. Such relationships were discovered in \cite{CN3}. 

 \begin{theorem}\label{T:5.1} Let $G$ and $G'$ are two compact Lie groups. If there exists a $k$-fold covering morphism $f:G'\to G$ for some odd $k$, then we have
\begin{equation} r_{2}G'=r_{2}G.\end{equation}  
\end{theorem}

As an application of Theorem \ref{T:5.1}, we have the following.

\begin{corollary} The $2$-rank $r_{2} G$ of $G$ depends only on the local class of $G$ if $G$ is one of the compact Lie groups $SU(k)$ and $E_{6}$, where $k$ is odd.
\end{corollary}

This corollary  follows from Theorem \ref{T:5.1} and the fact that the fundamental group of the bottom groups $E_{6}^{*}$ and $SU(k)^{*}$ have odd order with odd $k$.

For double covering Lie groups we have the following result.

 \begin{theorem}\label{T:5.2} Let $G$ and $G''$ are compact Lie groups. If $G$ is a double covering group of $G''$, then we have\begin{equation} \label{5.2} r_{2}G \leq 1+r_{2}G''.\end{equation} 
 \end{theorem}

This theorem was proved as follows: Let $A_{2}$ be a maximal antipodal set in $G$. The union $A_{2}\cup \gamma(A_{2})$ is also an antipodal set, where $\gamma$ is the covering transformation for the covering morphism $\pi: G\to G''$. Thus $\gamma$ stabilizes $A_{2}$. Consequently, $\pi(A_{2})$ is antipodal in $G''$. Hence we find
$$\#(A_{2})= 2\#(\pi A_{2})\leq  2\#_{2}G'',$$ which implies inequality \eqref{5.2}.

 \begin{remark} Inequality \eqref{5.2} is sharp, since the equality  holds in case $G$ is $SO(2m)$ for $m>2$.
 \end{remark}

\section{2-ranks for compact Lie groups with poles}

Suppose a finite group $\Gamma$ is acting on two symmetric spaces $M$ and $N$ freely as automorphism groups. Then $\Gamma$ acts on the product space
$M\times N$ freely. And the orbit space $(M\times N)/\Gamma$ is called the {\it dot product} of $M$ and $N$ (with respect to $\Gamma$) and denoted by $M\cdot N$. 

In most cases, $\Gamma$ will be the group of order two acting on $M$ and $N$ as the covering transformation groups for double covering morphisms. In the sequel, $\Gamma$ will not be mentioned in that case, if $\Gamma$ is obvious or if $\Gamma$ need not be specified.

\begin{example}\label{E:6.1} $SO(4)=S^{3}\cdot S^{3}$ and $U(n)=SO(2)\cdot SU(n)$. Here $\Gamma$ for $U(n)$ is the center of $SU(n)$, a cyclic group of order $n$.
\end{example}

The following two results from \cite{CN3} relate the 2-ranks, poles, and the dot product for compact Lie groups.

\begin{theorem}\label{T:6.1} If a compact connected Lie group $G$ has a pole, then 

\begin{itemize}
\item[{\rm (i)}] $\#_{2}(U(1)\cdot F)$ equals either $2^{r_{2} G}$ or $2^{1+r_{2} G}$;

\item[{\rm (ii)}] If the second case in {\rm (i)} occurs, then the centrosome of $G$ has the 2-number equal to $\#_{2}G$. 
\end{itemize}
\end{theorem}

\begin{theorem}\label{T:6.2} One has 
$$r_{2}G\leq r_{2}(Sp(1)\cdot G)\leq 2+r_{2}G$$
for a compact Lie group $G$ with a pole.
\end{theorem}

Theorem \ref{T:6.1} and Theorem \ref{T:6.2} were used in the determination of the 2-ranks for exceptional Lie groups

 \section{2-rank and Krull dimension}
 
 In commutative algebra, the {\it Krull dimension} of a ring $R$, named after Wolfgang Krull (1899-1971), is the supremum of the number of strict inclusions in a chain of prime ideals, not the number of primes. More precisely, we say that a strict chain of inclusions of prime ideals of the form: 
 $${\mathfrak p}_{0}\subsetneq {\mathfrak p}_{1}\subsetneq \cdots \subsetneq {\mathfrak p}_{n}$$
 is of length $n$. That is, it is counting the number of strict inclusions. 
   
Given a prime ideal ${\mathfrak p}\subset R$, we define the {\it height} of ${\mathfrak p}$ to be  the supremum of the set 
$$\{n\in {\mathbb N}: {\mathfrak p} \hbox{ is the supremum of a strict chain of length $n$}\}.$$
Then the Krull dimension is the supremum of the heights of all of its primes.
 
 Let $G$ be a compact Lie group. Put $$H_{G}^{*}=H^{*}(BG; {\mathbb Z}_{2}),$$ where $BG$ is a classifying space for $G$. Let $N_{G}^{*}\subset H_{G}^{*}$ denote the ideal of nilpotent elements.  Then $H_{G}^{*}/N_{G}^{*}=H_{G}^{\#}$  is a finitely generated commutative algebra.  
 
 In \cite{Q1}, D. Quillen studied the relationship between $H_{G}^{\#}$  and the  structure of the Lie group $G$. In particular, he proved that, under some suitable assumptions, the Krull dimension of $H_{G}^{\#}$ is equal to the 2-rank of $G$. He proved the result by calculating the mod 2 cohomology ring of extra special 2-groups. 
 
 Quillen's result provides an affirmative answer to a conjecture of M. F.  Atiyah and R. G. Swan (cf. \cite{At,Sw}).

\section{2-ranks of dihedral groups}

We denote by $D[4]$ the dihedral group of order 8, or the
automorphism group of a square in the plane. Thus $D[4]$ is generated by the
reflections in the $x$-axis and the line $y=x$ in the Euclidean plane. Clearly \begin{equation}\#_{2}D[4] = 4.\end{equation}

Let $Q[8]$ denote the quaternion group, generated, by $i$ and $j$ in the group of the
nonzero quaternions, where $i$ and $j$ together with $k$ form a standard basis for the
pure quaternions. One has \begin{equation}\#_{2}Q[8] = 2.\end{equation} Their commutator subgroups have the
2-numbers both equal to 4.

\section{2-ranks of classical groups} 

In this section, we present the 2-ranks of all classical Lie groups from \cite{CN3}. For the simply-connected classical Lie groups $U(n),\, SO(n)$ and $Sp(n)$, this was already done by Borel and Serre in \cite{BS}.

 \begin{theorem}\label{T:9.1} Let $U(n)/{\mathbb Z}\mu$  by the quotient group of the unitary group $U(n)$ by the cyclic normal subgroup ${\mathbb Z}\mu$ of order $\mu$. Then we have 
  \begin{equation} r_{2}(U(n)/{\mathbb Z}\mu)= 
 \begin{cases} n+1 & \text{if $\mu$ is even and $n=2$ or $4$;}\\ n & \text{otherwise.}\end{cases} \end{equation}
 \end{theorem}

The proof of this theorem based on several rather complicate lemmas whose proofs used linear algebra.

 \begin{theorem}\label{T:9.2} For $SU(n)/{\mathbb Z}\mu$, we have 
 \begin{equation} r_{2}(SU(n)/{\mathbb Z}\mu)= 
 \begin{cases} n+1 & \text{for $(n,\mu)=(4,2)$;}
 \\ n  & \text{for $(n,\mu)=(2,2)$ or $(4,4)$;}
 \\ n-1 & \text{for the other cases.} \end{cases} \end{equation}
 \end{theorem}

The proof of Theorem \ref{T:9.2} is similar to that of Theorem \ref{T:9.1}.

The next three results are concerned with other classical groups and their adjoint groups. Basically, they are proved with the same method as Theorem \ref{T:9.1}.

 \begin{theorem}\label{T:9.3} One has $r_{2}(SO(n))=n-1$ and, for $SO(n)^{*}$, we have
 \begin{equation} r_{2}(SO(n)^{*})= 
 \begin{cases} 4 & \text{for $n=4$;}\\ n-2 & \text{for $n$ even $>4$}.\end{cases} \end{equation}
 \end{theorem}

\begin{remark} $SO(n)^{*}$ is $SO(n)$ for $n$ odd and $SO(n)/\{\pm 1\}$ if $n=2n'$ is even $>2$. 
\end{remark}

 \begin{theorem}\label{T:9.4} Let $O(n)^{*}=O(n)/\{\pm 1\}$. We have
 
 {\rm (a)} $r_{2}(O(n))= n$;
 
 {\rm (b)} $r_{2}(O(n)^{*})$ is $n$ if $n$ is 2 or 4, while it is $n-1$ otherwise.
 \end{theorem}

For $Sp(n)$ and $Sp(n)^{*}$ we have the following.

 \begin{theorem}\label{T:9.5} One has $r_{2}(Sp(n))=n$, and, for $Sp(n)^{*}$, we have
  \begin{equation} r_{2}(Sp(n)^{*})= \begin{cases} n+2 & \text{for $n=2$ or $4$}\\ n+1 & \text{otherwise}.
\end{cases} \end{equation}
Thus we also have
 \begin{equation}r_{2}(Sp(n)^{*})= r_{2}(U(n)/{\mathbb Z}_{2}) + 1 \end{equation} for every $n$. 
\end{theorem}

\section{2-ranks of spinors and semi-sipnors}

Now, we consider the spinor $Spin(n)$ and related groups. Recall that  $Spin(n)$ is a subset of the Clifford algebra $Cl(n)$, which is generated over $\mathbb R$ by
the vectors ${\bf e}_{i}$ in the fixed orthonormal basis of $\mathbb R^{n}$ and subject to the conditions
${\bf e}_{i}{\bf e}_{j}=-{\bf e}_{j}{\bf e}_{i}$, and ${\bf e}_{i}{\bf e}_{i}= -1,\, i\ne j$ (cf. \cite{ABS,C}). 

Under the projection $$\pi:Spin(n)\to SO(n),$$ the member $(\cos\theta) {\bf e}_{1}+(\sin\theta){\bf e}_{2}$ of $Spin(n)$, for instance, projects to the rotation of the $({\bf e}_{1},{\bf e}_{2})$-plane by the angle $2\theta$, carrying ${\bf e}_{1}$ into 
$(\cos 2\theta){\bf e}_{1}+(\sin 2\theta){\bf e}_{2}$, for every real number $\theta$, where $\pi$ denotes the double covering homomorphism with the kernel $\{1,-1\}$. 

Every maximal antipodal group $A_{2}$ in the spinor $Spin(n)$ projects into a diagonalizable subgroup of $SO(n)$. Therefore, we may assume that $A_{2}$ is
a subgroup of $$E(n)= \{\pm {\bf e}_{I}: I\subset \{1,2,\ldots,n\}\},$$ where ${\bf e}_{I} = {\bf e}_{i}{\bf e}_{j}\cdots {\bf e}_{k}$ ($=1$ if $I$ is empty) for any subset $$I = \{i,j, \ldots,k\},\;\;\; i<j<\cdots <k$$ of $\{1,2,\ldots,n\}$. 

For $Spin(n)$ we have the following results from \cite{CN3}.

 \begin{theorem}\label{T:10.1} We have
 \begin{equation}\notag
  r_{2}(Spin(n)) =\begin{cases} r+1 & \text{if  $n \equiv -1,0$ or {\rm 1 (mod 8)}}\\ 
 r &\text{otherwise},\end{cases}
 \end{equation} where $r$ is the rank of $Spin(n)$, $r = [n/2]$. 
 \end{theorem}

 \begin{theorem}\label{T:10.2} {\rm (PERIODICITY)} One has
$$ r_{2}(Spin(n+8))= r_{2}(Spin(n))+4$$
for $n\geq 0$. 
\end{theorem}

The group $Pin(n)$ was introduced by M. F. Atiyah, R. Bott and A. Shapiro in \cite{ABS} while they studied Clifford modules. $Pin(n)$ is a group in the Clifford algebra $Cl(n)$ and it double covers $O(n)$ and whose connected component $Spin(n)$ double covers $SO(n)$.

 \begin{theorem}\label{T:10.3} For $Pin(n)$ we have $$r_{2}(Pin(n))= r_{2}(Spin(n + 1))$$ for $n\geq 0$.
 \end{theorem}

The proofs of Theorem \ref{T:10.1}, Theorem \ref{T:10.2} and Theorem \ref{T:10.3} base mainly on the following lemma.

\begin{lemma}\label{L:10.1} We have the following:

\begin{itemize}
\item[(a)] One has ${\bf e}_{I}{\bf e}_{J}={\bf e}_{J}{\bf e}_{I}$ if and only if the cardinalities satisfy
$$(\#I)(\#J)-\#(I\cap J)\equiv 0\;\;  ({\rm mod} \,2);$$

\item[(b)] one has $({\bf e}_{I})^{2}=1$ if and only if $\#I\equiv 0$ or
$3$ $({\rm mod}\, 4)$; and 

\item[(c)] ${\bf e}_{I}$ is a member of $Spin(n)$ if and only if $\#I$ is even.
\end{itemize}
\end{lemma}

Let $SO(4m)^{\#}$ denote $Spin(4m)/\{1, e_{((4m))}\}$ the semi-spinor group.

For the semi-spinor group $SO(4m)^{\#}$, we have the following result from \cite{CN3}.

 \begin{theorem}\label{T:10.4}   We have
 \begin{equation}\notag
  r_{2}(SO(4m)^{\#}) =\begin{cases} 3 & \text{if  $m=1$}\\  6 &\text{if $m=2$,}\\ r+1 &\text{if $m$ is even $>2$, }\\
 r &\text{if $m$ is odd $>1$},\end{cases}
 \end{equation}
 where $r$ is the rank $2m$ of $SO(4m)^{\#}$.  
 \end{theorem}

For the proof of Theorem \ref{T:10.4}, we have applied the dihedral groups $D[4]$ and $Q[8]$ (see \S 8) and we also investigated the inverse image of a maximal antipodal subgroup $A_{2}$ of $SO(4m)^{\#}$ under the projection: $Spin(4m)\to SO(4m)^{\#}$.

\begin{remark} The maximal 2-subgroups of both $Spin(16)$ and  $SO(16)^{\#}$ have been studied independently by J. F. Adams in \cite{Ad}.
\end{remark}

\section{2-ranks of Exceptional groups} 

Finally, we provide the 2-ranks of compact exceptional Lie groups from \cite{CN3}. For $G_{2}$ and $F_{4}$ this was already done in \cite{BS}.

\begin{theorem}\label{T:11.1} One has $$r_{2}G_{2}=3,\;\;  r_{2} F_{4}=5,\;\;  r_{2}E_{6}=6,\;\; r_{2}E_{7}=7
,\;\; r_{2}E_{8}=9$$ for the simply-connected exceptional simple Lie groups.\end{theorem}

For the proof of Theorem \ref{T:11.1}, we have applied Theorem \ref{T:6.2} and the $(M_{+},M_{-})$-theory of compact symmetric spaces from \cite{CN1}, as well as linear algebra.

For the bottom space $E^{*}_{6}$, we have the following result from \cite{CN3}.

 \begin{theorem}\label{T:11.2} One has $r_{2}E^{*}_{6}=6$. \end{theorem}

This theorem follows from Corollary \ref{T:6.1} and  Theorem \ref{T:11.1}.

\begin{remark} Independently, J. F. Adams studied  in \cite{Ad} the 2-subgroups of $E_{8}$. He proved that a maximal 2-subgroup of the compact exceptional simple group $E_{8}$ fall into just two conjugacy classes $D(T^{8})\subset Spin(16)$ and $EC^{9}$ (see \cite{Ad} for details).

\end{remark}
\section{Concluding remarks}

 (1) For every complex flag manifold $M_{\mathbb C}$ there exists a positive integer $k_{0}=k_{o}(M_{\mathbb C})\geq 2$ such that for each integer $k\geq k_{0}$ there exists a $k$-symmetric structure on $M_{\mathbb C}$. 
 
To extend the notion of 2-number $\#_{2}M$ of Chen-Nagano for compact symmetric spaces,  C. U. S\'anchez introduced in \cite{Sa} the $k$-number $\#_{k}(M_{\mathbb C})$ of the complex flag manifold $M_{\mathbb C}$ as the maximal possible cardinality of the so-called $k$-sets $A_{k}\subset M_{\mathbb C}$ with the property that for each $x\in A_{k}$ the corresponding $k$-symmetry fixes every point of $A_{k}$. He then proved that
$$ \#_{k}(M_{\mathbb C})=\dim H^{*}(M_{\mathbb C},\mathbb Z_{2})$$
for each complex flag manifold $M_{\mathbb C}$. 

Using the fact that each real flag manifold $M$ can be isometrically embedded into a complex flag manifold $M_{\mathbb C}$  (the so-called complexification of $M$), he defined in \cite{Sa2}  the index $p$ of $M$ as the smallest prime number $p\geq k_{0}(M_{\mathbb C})$. Moreover, he defined the {\it index number} $\#_{I}M$ of $M$ as the maximal possible cardinality of the $p$-sets $A_{p}M$. 

C. U. S\'anchez proved that
$$\#_{I}M=\dim H^{*}(M,\mathbb Z_{2})$$
for any real flag manifold $M$.
 
(2)  B.-Y. Chen conjectured in \cite{c} that the 2-number $\#_{2}M$ of a compact symmetric space (or a compact Lie group) $M$ is closely related to the smallest number of cells that are needed for a $CW$-complex structure on $M$.

  By direct computation one can see that the 2-number $\#_{2}M$ is in fact equal to the smallest number of cells needed to have a $CW$-complex structure of $M$ if $M$ is a sphere, a real projective space or a Hermitian symmetric space.
 
J. Berndt, S. Console and A. Fino proved in \cite{BCF} that the index number $\#_{I}M$ of a real flag manifold $M$ coincides with the smallest number of cells that are needed for a $CW$-complex structure on $M$. Furthermore, they showed that the intersection of the fixed point sets of all $k$-symmetries for all $k\geq k_{0}$ at a point $x$ of a complex flag manifold coincides with the critical point set of the height function with respect to $x$.
 Moreover, they proved  that the index number $\#_{I}M$ is the number of cells of the Bruhat decomposition of $M$ which is determined by the above generic height functions.
 
 (3) H. Tasaki  investigated in \cite{Ta} the 2-number of the complex hyperquadric 
 $$Q_{n}(\mathbb C)=\{[z_{1},z_{2},\ldots,z_{n+2}]\in {\mathbb C}P^{n+1}: z_{1}^{2}+z_{2}^{2}+\cdots +z_{n+2}^{2}=0\}.$$ It is known that $Q_{n}(\mathbb C)$ is holomorphically isometric to the Hermitian symmetric space $SO(n+2)/SO(2)\times SO(n)$, which is the Grassmann manifold $\tilde G_{2}({\mathbb R}^{n+2})$ consisting of all oriented linear subspaces of dimension 2 in ${\mathbb R}^{n+2}$.
This Hermitian symmetric space  admits certain real forms $S^{k,n-k}$ defined by $$S^{k,n-k}=(S^{k}\times S^{n-k})/{\mathbb Z}_{2}.$$  Let $k,\ell$ be integers with $0\leq k\leq \ell \leq [n/2]$, and let $L_{1}$ and $L_{2}$ be real forms of $\tilde G_{2}({\mathbb R}^{n+2})$ which are congruent to $S^{k,n-k}$ and $S^{\ell,n-\ell}$, respectively. 
In \cite{Ta}, H. Tasaki proved that if $L_{1}$ and $L_{2}$ intersect transversally, then $L_{1}\cap L_{2}$ is a maximal antipodal set of $L_{1}$ and an antipodal set of $L_{2}$. Moreover, if $k=\ell=[n/2]$,  then $L_{1}\cap L_{2}$ is a maximal antipodal set of $\tilde G_{2}({\mathbb R}^{n+2})$. As a consequence, he proved that any real form of $\tilde G_{2}({\mathbb R}^{n+2})$ is a globally tight Lagrangian submanifold.

(4) By computing the geometric invariant $\#_{2}G$ of a compact Lie group $G$, we are able to determine the 2-rank $r_{2}G$ of all compact simple Lie groups $G$ via Theorem \ref{T:3.1} and applying the $(M_{+},M_{-})$-method. 
It would be quite interesting to establish a general method for computing the $p$-rank $r_{p}G$ of a compact Lie group $G$ for a prime number $p>2$ in the spirit of Theorem \ref{T:3.1}.

\end{document}